\documentclass[11pt,a4paper]{article}

\usepackage{amsmath,amssymb}
\usepackage{verbatim}
\usepackage{enumerate}
\usepackage{hyperref}
\usepackage{fontenc} 
\usepackage[utf8]{inputenc}
\usepackage{authblk}
\usepackage{a4,graphicx}
\usepackage{float}
\usepackage{caption}
\usepackage{subcaption}

\title{Math That Matters:\\
\large Enhancing Academic  Mathematics'  Impact on Society\\
\large \textit{Plenary address to 2019 ICAPTA Conference, Abuja Nigeria} }
\date{\vspace{-5ex}}
\author[1]{\small Christopher Thron}
\author[2]{\small Monira Taj Elsir Hamid Ali }

\affil[1]{\small Texas A\&M University-Central Texas, Killeen, TX 76548, USA}
\affil[2]{\small Alzaiem Alazhari University, Khartoum Sudan}

\affil[1]{\small email address: thron@tamuct.edu}
\affil[2]{\small email address: monirahamid@gmail.com}

\begin{document}

\maketitle

\noindent{\bf Abstract:}  
Most if not all of today’s revolutionary technologies have a common foundation, namely the intelligent use of information. It is clear that computers play a central role: but the contribution of mathematics, though less visible, is no less critical. The conceptual tools and insights provided by mathematics are the keys to unlocking the full information-processing potential of computers. Academic mathematicians now have an unparalleled opportunity to make a huge impact on modern society: but to take advantage of this opportunity, mathematicians must prioritize making these tools and concepts accessible to a wider audience.

In this paper we present three examples of mathematics with significant social benefit: Dmitri Bertsimas’ study of diabetes using $k$-nearest-neighbor methodology;  Development of mathematical software (Matlab and Sage); and ongoing development of data representation and visualization software to facilitate analysis of survey data.  We also suggest steps to be taken by academic mathematicians  in Nigeria towards enhancing the positive impact of mathematics on society.  
\\



\section{Introduction}

\subsection{A Richer Picture of Academic Mathematics}
What do mathematicians do? Ask a pure mathematician, and he may answer: “Prove theorems”. Ask an applied mathematician, and he may say, “Find solutions to mathematical models, and/or investigate their properties.”

\begin{figure}[H]
	\centering
	\includegraphics[width=7in]{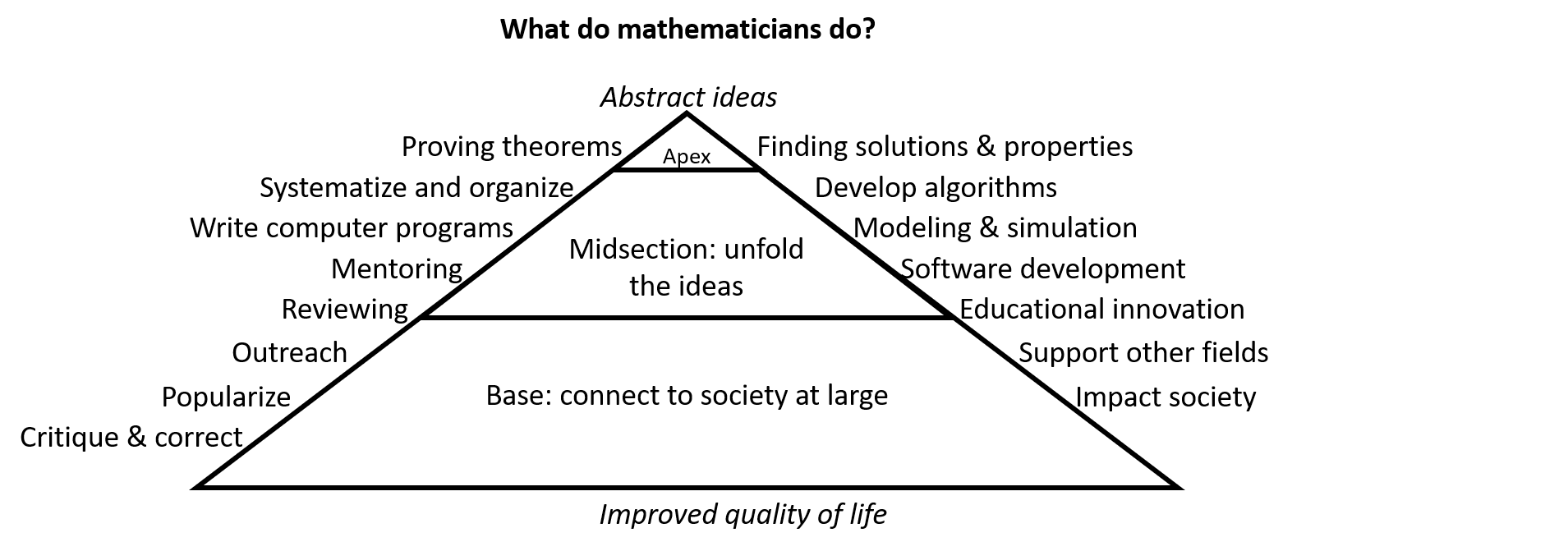}
	\caption{A ``mathematical pyramid'', showing a range of constructive activities in academic mathematics.}
	\label{fig:pyramid}
\end{figure}

In this paper we would like to suggest a somewhat richer picture of the academic mathematician’s job description. 
Figure~\ref{fig:pyramid}  is a representation of our so-called “mathematical pyramid”. At the apex are the conventional activities attributed to mathematicians: proving theorems and finding solutions. But beneath the apex, many other activities are included. The pyramid resembles the organizational structure of high-tech companies: a preponderance of secondary activity is required to make the brilliant ideas at the top connect with society to bring the greatest possible benefit. 

Should mathematicians be concerned about impact? G.H. Hardy was the most famous British mathematician in the early 20th century, and  has been made famous to the modern public by the film ``The Man Who Knew Infinity''. G. H. Hardy was famous for saying mathematics should be entirely impractical.  Here are some quotes from Hardy's book, \emph{A Mathematician's Apology}\cite{Hardy}: 
\begin{quote}
I have never done anything ``useful''. No discovery of mine has made, or is likely to make, directly or indirectly, for good or ill, the least difference to the amenity of the world.

No one has yet discovered any warlike purpose to be served by the theory of numbers or relativity, and it seems unlikely that anyone will do so for many years.

We have concluded that the trivial mathematics is, on the whole, useful, and that the real mathematics, on the whole, is not.
\end{quote}
Over time, Hardy has been proved wrong—not just a bit off, but totally mistaken. The theory of numbers is the foundation for cryptography, which is the basis for modern internet commerce; and Einstein's relativistic mass equation $E=mc^2$ expresses the mass-energy equivalence that powers the atomic bomb. More constructively, the theory of general relativity is essential to the proper functioning of the Global Positioning System (GPS). By general consensus, the most beautiful and mysterious of mathematical equations is 
\begin{equation}
e^{i\pi} = -1,
\end{equation}
which is a simple corollary of the properties of complex exponentials that are essential tools in modern communications. As to Hardy's own work, Hardy-Weinberg equilibrium is an important result in modern genetics (although to be fair this does not rely on the ``interesting'' advanced mathematics that was Hardy's main focus). 

Hardy’s opinion may be viewed as reflection of his atheistic worldview. He believed that the imagination of man surpassed the creativity of the Creator. In his own words: `` `Imaginary' universes are so much more beautiful than this stupidly constructed ‘real’ one.”''

In stark contrast to Hardy, most if not all of the greatest mathematicians have been concerned with practical impact. Isaac Newton did not even consider himself as a mathematician but as a “natural philosopher”. At the age of  56 he became Master of the Royal Mint, and that capacity he stamped out counterfeiting and put English currency on a solid basis, paving the way for economic growth and stability.   C. F. Gauss was actually a professor of astronomy, not of mathematics. Over several years he conducted a geodetic survey of the Kingdom of Hanover, traveling many miles by horseback--and out of this project came the subject of differential geometry. John Von Neumann made foundational contributions to economics, and has been called the father of the modern computer.  Professor Terry Tao of  U.C. Berkeley is a strong contender for the title of “greatest mathematician alive”. His areas of interest according to Wikipedia are “harmonic analysis, partial differential equations, combinatorics, compressed sensing and analytic number theory.”  His work in compressed sensing has spawned a whole new technology, including what is called the “1-pixel camera”.\cite{Tao} 

A number of African mathematicians have received support from the Simons foundation, which gives generous research grants to mathematicians worldwide\cite{Simons0}. The foundation was started and funded by James Harris Simons, a former math professor (and co-discoverer of Chern-Simons forms in algebraic topology) who went on to become a financial trader on Wall Street.\cite{Simons}  

In summary it seems clear that, although not all mathematicians need to be conducting impactful research, we certainly need to encourage those that do. Their work can bring enormous benefit to the entire field, including those who inhabit the apex of the ``pyramid''. If mathematics is widely seen to be beneficial, funding will increase; but if it makes little connection with society, funding will be neglected.

Let me mention briefly another mathematician who has performed great services for mathematics, apart from proving theorems and solving problems. Professor Steven Strogatz of Cornell University is a leading authority in nonlinear dynamics and complex systems. He is co-author of a paper entitled “Collective dynamics of small-world networks”, which so far has more than 38,000 citations in various scientific papers. He has also written five books: ``Nonlinear Dynamics and Chaos'', ``Sync: How Order Emerges From Chaos in the Universe, Nature, and Daily Life'', ``The Calculus of Friendship'', ``The Joy of x'',  and ``Infinite Powers''.  As perhaps you can tell, the last four books have targeted non-mathematician audiences.  He also has written several articles on mathematics for the New York Times, which is one of the most prestigious newspapers in the U.S. and is read nationwide.  In short, professor Strogatz is a math popularizer. 

Personal experience leads me (C.T.)  to believe that this kind of popularization work is sorely needed in Nigeria. In December 2018 I   had the great privilege of participating in the first ``International Conference on Mathematical Modeling of Environmental Pollution'' at the National Mathematical Center in Abuja . During the conference I was interviewed by a Nigerian reporter who asked me, ``What  can mathematics do to help the environment?'' She was clearly  baffled that mathematics could have anything at all to do with reducing pollution. To be honest, I was not prepared for the question and didn't give a very good answer. I encourage all of us to prepare ourselves: in another context that some readers may be familiar with, we are exhorted to ``be ready always to give an answer to every man that asketh you a reason of the hope that is in you''\cite{Peter}. A similar exhortation applies to our mathematical occupation.

 I expect that there are mathematicians who are already doing popularization work in Nigeria. Their work should be recognized and supported by the entire academic mathematical community.  I  hope that more will step up and do the same. I encourage readers to look up on the Internet the Youtube channel ``Numberphile''\cite{Numberphile}, which may give you some ideas of the possibilities.

\section{Examples of High-Impact Mathematics}

Let us now turn to consider three different examples of high-impact mathematical research. 

\subsection{Mathematics and Treatment of Diabetes}
Our first example of impactful mathematics is work done by Professor Dmitri Bertsimas and colleagues at MIT. He was moved to study diabetes because the disease runs in his own family: his mother and aunt were both diabetic. Bertsimas was impressed by the fact that although they had the same disease, their prescribed treatments were very different.  He became interested in finding a method to identify the very best treatment for each individual patient. So he arranged with Boston Medical Center to obtain electronic records from 11,000 anonymous diabetic patients. The records included the patients’ demographic and medical history data. Bertsimas and his colleagues distinguished 13 different therapies administered to different patients at various times. Each patient was represented by a point in a multidimensional  space, where each dimension in the space represented a different demographic or medical history characteristic of the patient. A metric (i.e. a notion of distance) was defined on this space, so that the closer the distance between two patients, the greater the overall similarity between their characteristics. 

In order to predict the effect of the 13 different therapies on a particular patient, for each therapy the 30 therapy recipients who were closest to the patient in question were identified. A weighted average of the outcomes of these thirty patients was taken as a best indicator for the effectiveness of the corresponding therapy. (Readers familiar with machine learning may recognize this as a variation of the standard classification technique known as “$k$-nearest neighbor”.) As a result, a predicted best treatment tailored for the individual patient could be identified. Professor Bertsimas  reported that a  significant lowering of hemoglobin A1c levels (an indicator of blood glucose levels) was observed when patients were put on the predicted best treatment\cite{Bertsimas}. 

\begin{figure}[H]
	\centering
	\framebox{\includegraphics[width=4in]{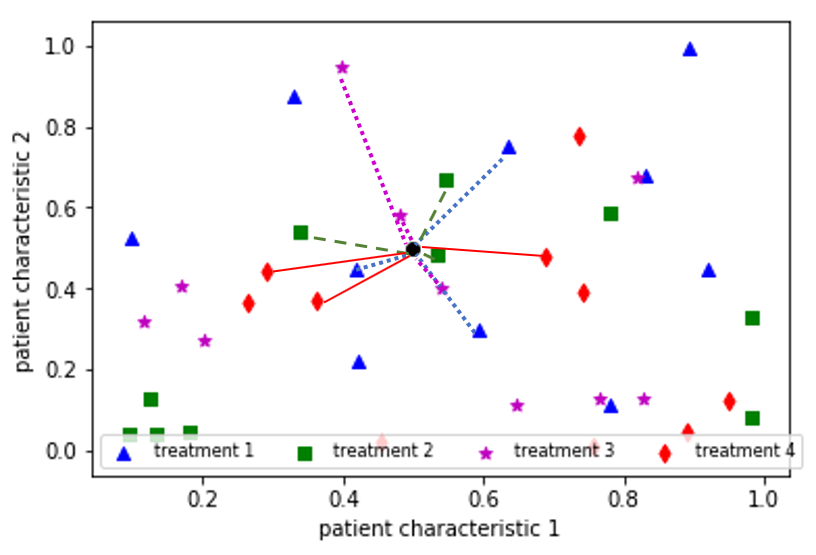}}
	\caption{Schematic illustration of Bertimas'  ``nearest neighbor'' method, with 3 nearest neighbors taken for each therapy (Bertsimas actually used 30).}
	\label{fig:knn3}
\end{figure}

A schematic illustration of Bertsimas' method is shown in Figure~\ref{fig:knn3}.  In the figure, the round black point signifies the "patient" of interest.  The three nearest neighbors for four different therapies are recognized, which will yield four different effectiveness estimates based on weighted averages of the neighbors' outcomes. The therapy that is the most effective of the four is the recommended therapy for the patient.  

This is just one example of mathematics applied to the health field.  This is a burgeoning area of research, because of the availability of new techniques to analyze massive amounts of data. I would encourage any researcher that if you find an opportunity to jump into this field, by all means do so. Right now is a very exciting time in this new area, and there is a prime chance to make contributions of fundamental importance.

\subsection{Development of Mathematical Software}
\begin{figure}[H]
	\centering
	\framebox{\includegraphics[width=3in]{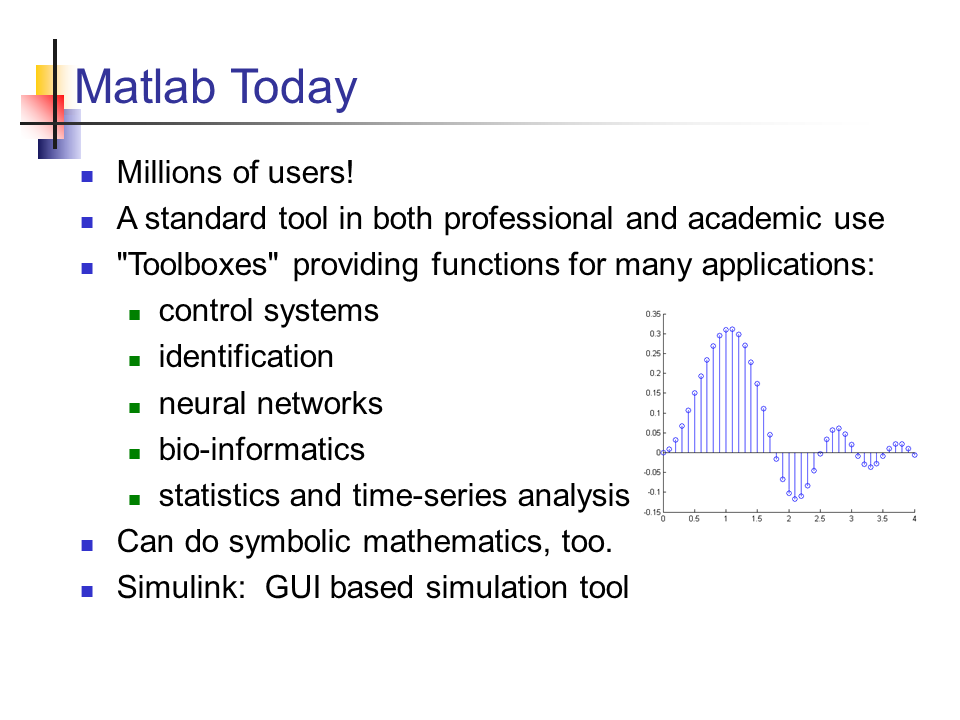}}
	\caption{Summary of Matlab's impact on modern science and technology (from \cite{Brucker}).}
	\label{fig:matlab}
\end{figure}

The second example is actually two examples in one, and both involve the development of software. Most readers will have heard of MATLAB. Even those who have never used MATLAB, have been most likely been impacted by it if they have done any programming,  because other software has copied its features and capabilities (for example, the syntax of Octave and the plotting routines in Python were imitated from MATLAB). Figure~\ref{fig:matlab} summarizes the main features of MATLAB.

Matlab originated in the work of Cleve Moler, who at the time was professor of mathematics at the University of New Mexico. Before developing MATLAB he had already developed LINPACK, which is a Fortran library of numerical linear algebra routines (and by the way, was supported by grants from the  U.S. government's National Science Foundation). In order to increase the useability for his students, he developed MATLAB (Matrix Laboratory) as a simple, interactive interface to access the powerful numerical routines in LINPACK.   Apparently the mathematics and computer science students in his classes weren't  impressed, but the engineering students were enthusiastic. A control engineer named Jack Little saw the capabilities of MATLAB, and approached Moler about rewriting it to be faster and more useable. Together with Steve Bangert they incorporated MathWorks in 1984, and the company grew steadily from there.\cite{matlabHistory}  In 1992, MathWorks had 100 employees, and today it has over 1000\cite{matlabCo}. Moler left academia in 1989 to become chief mathematician at MathWorks. Financially this was not a bad move for Moler--according to the Bloomberg website, he received a salary plus bonus of over \$500,000 last year\cite{Bloomberg}.

Note that MATLAB originated in academic research. Through partnership with non-academics, it became commercialized. Eventually the development became too big, and Moler had to leave academics and transition to industry. This is a very healthy trend. It is quite common in academic biomedical and electrical engineering departments in the U.S., and  in other tech-savvy countries such as India\cite{sciencemag,ecoTimes}. A free interchange between industrial and academic mathematics and mathematicians would breathe new life into academic mathematics in Nigeria. Computational software is one promising area where this interchange could take place.

Besides Matlab, another popular mathematical software package is Sage. Sage actually has more abstract algebraic capabilities than Matlab.  With Sage you can compute permutation groups, define p-adic numbers, algebraic varieties and elliptic curves over rings, Modular forms, Dirichlet characters, and so on. Sage pulls together a large number of existing specialized mathematical softwares including GAP for computer algebra, PARI for elliptic curves, and so on\cite{Sage}. 

\begin{figure}[H]
	\centering
	\includegraphics[width=3.5in]{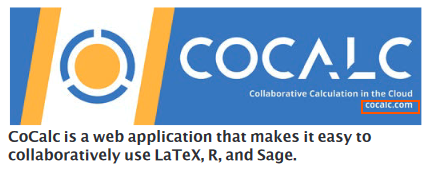}
	\caption{Logo for CoCalc, the main offering of SageMath Inc (from \cite{Stein}).}
	\label{fig:coCalcLogo}
\end{figure}

The primary developer of Sage is William Stein, formerly professor of mathematics at the University of Washington. His experience in developing Sage as an academic was somewhat rockier than Moler's. Stein's  vision was (and is) to develop an open-source alternative to Matlab and other proprietary mathematical software systems. He originally planned to accomplish this development as an academic, and assembled a talented group of Ph.D. students to work on Sage. Unfortunately, once all his students graduated they all left and got high-paying jobs at Google, Facebook, Microsoft, or similar companies. He also encountered significant resistance within the academic community\cite{Wing}. In Stein's own words\cite{Stein}: ``Some of the mathematical community genuinely considers the creation of mathematical software as not being `real mathematics'. If you most love improving mathematical software, you will probably end up having to leave academia (pure math, at least).''  Nonetheless, Sage has made a huge impact within pure mathematics: it has been cited in 357 articles,  39 theses, and 40 books\cite{Sage}.  Sage has also been integrated into online math textbooks in Linear Algebra and Abstract Algebra\cite{Beezer,Judson}. in 2016 Stein himself left academia and is now CEO of his own company. SageMath Inc. The company's main product is CoCalc, a browser-based  software for online mathematical collaboration (see Figure~\ref{fig:coCalcLogo}).

Here we mention briefly also Dr. Allan Steel at the University of Sydney.  Dr. Steel is not a professor, but rather a Senior Research Fellow and member of the Computational Algebra Research Group (I am not aware of any comparable type positions in Nigeria.) He is the principal developer of Magma, which is described on the Magma web page as follows: ``Magma is a large, well-supported software package designed for computations in algebra, number theory, algebraic geometry and algebraic combinatorics. It provides a mathematically rigorous environment for defining and working with structures such as groups, rings, fields, modules, algebras, schemes, curves, graphs, designs, codes and many others.  Magma also supports a number of databases designed to aid computational research in those areas of mathematics which are algebraic in nature.''\cite{Magma}  Magma had already been cited in about 7000 scientific publications.  Magma is not free: however, since 2013 the Simons Foundation (mentioned above) has covered the cost for U.S. educational institutions.   

Moler, Stein, and Steel should not be statistical outliers.  Mathematicians should by rights be in the forefront of this kind of software development. Mathematicians are uniquely able to lay down the theoretical foundations for efficient software. They may then partner with computer scientists, who can efficiently implement the mathematicians' conceptions. A recent example of this is the computer language Julia, which may eventually  replace Python and Matlab as the leading software for computational mathematics\cite{Julia1,Julia2} (for one, Julia executes much faster than Matlab or Python). Alan Edelman, another math professor at MIT and authority in random matrix theory, was one of the co-creators of Julia.      

\subsection{A Statistical `My Dear Watson' }
As a third and final example of (potentially) high-impact mathematics,  I will briefly present a current research project that is joint work with  Monira Taj Elsir Hamid Ali ,  a lecturer and Ph.D. candidate at Alzaeim Alazhari University in Khartoum, Sudan. First, let me present the motivation.

\begin{figure}[H]
	\centering
	\framebox{\includegraphics[width=3.5in]{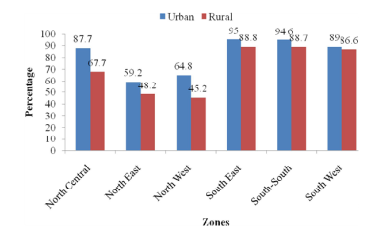}}
	\caption{Bar chart from NARHS report (2012)\cite{NARHS}. Note that the x-axis categories is arranged in alphabetical order, which is irrelevant to the data. A much better choice would be to order the categories in decreasing order.}
	\label{fig:NARHS}
\end{figure}

In the modern world, statistics is an essential tool in making intelligent policy decisions on the personal, corporate, and governmental levels. Unfortunately, many key decision-makers do not have the statistical knowledge required to make good use of important data they possess. They need the data to be represented in a way that they can understand.  The data science community in the U.S. has recognized this, and data visualization is now a large and active research area. It seems that this is an area in which Nigeria is woefully lagging. 
 For example, Nigeria has conducted three National AIDS and Reproductive Health Surveys (NARHS) in 2003, 2005, and 2012.  The 2012 report has 527 pages, and 17 figures\cite{NARHS}. All but one of the figures  was a simple bar chart (an example is shown in Figure~\ref{fig:NARHS}), and  the lone exception was a simple line chart. Most of the figures concerned general demographics of respondents.   Only two of them related  HIV infection to other factors: the factors were age, sex and geographical location. All other information was contained in tables full of incomprehensible numbers. It seems doubtful that the wealth of information in this bulky report is readily comprehensible by government officials and legislators who are be responsible for funding and determining public health policy.

Decision-makers who have sufficient funding may hire statisticians to analyze and present their data. However, the statisticians often do not understand the practical situation being analyzed, and so they may bring out aspects of the data which are irrelevant as far as policy is concerned.  There is a big difference between statistical significance and practical significance, and statisticians are prone to focusing on the former, without properly appreciating the latter.  

At this point we may explain the title of this section, ``A Statistical `My Dear Watson'''. Sherlock Holmes is a fictional character (created by Arthur Conan Doyle) known as a brilliant detective who is able to make incredibly detailed, correct deductions from very small clues. In the Sherlock Holmes books, Doctor Watson is Holmes' friend,
who (although much less brilliant that Holmes)  has a unique ability to stimulate Holmes' ideas through his comments. When Watson's perplexities succeeded in giving Holmes a new idea, Holmes would exclaim, "Elementary, my dear Watson''. Hence, the name for our project.  

Our idea is to create  software that functions as Doctor Watson did: namely to stimulate users (and especially users with little mathematical background) to get ideas and draw conclusions from data, based on easily-understood graphical representations. The steps to our approach are as follows. 
\begin{enumerate}
\item
Generate a library of plots (as .pdf files)  that represent different aspects of the data, and in particular relationships between the different variables;
\item
Supply the user with a set of leading questions which guides the user's understanding of the significance of the data;
\item
Enable the user to perform simple data manipulations (like removing, adding, or merging specific categories of the different variables) and generate other plots to answer questions raised by the original library of plots. 
\end{enumerate}

\begin{figure}[H]
	\centering
	\framebox{\includegraphics[width=4.5in]{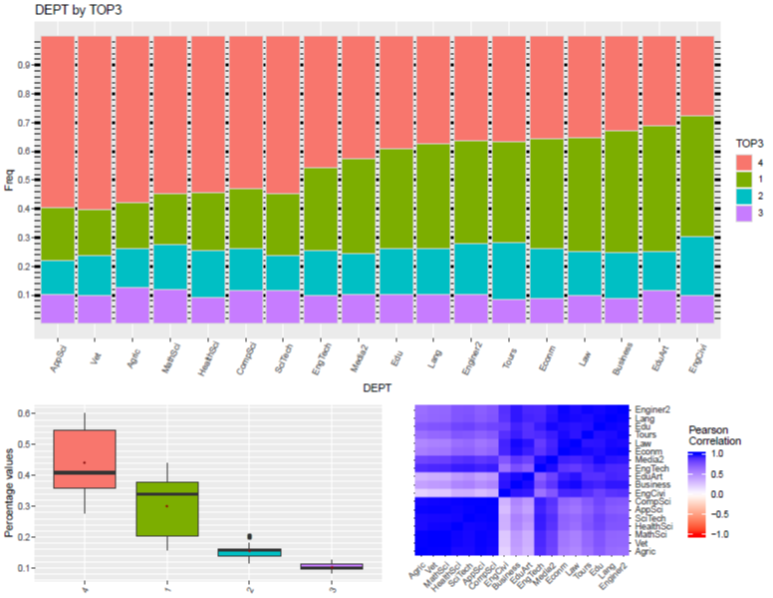}}
	\caption{Two-variable plot produced by My Dear Watson, based on Sudan university admissions data (2015-16).  The bar variable is university department, and the color variable indicates whether the admitted student obtained his/her first, second, third, or higher choice of university and department.}
	\label{fig:2varPlot}
\end{figure}

Figure \ref{fig:2varPlot} gives an example of a 2-variable display generated by My Dear Watson (which is written in R, and is accessed through the RStudio graphical interface). The display includes a percentage bar plot, a box plot, and a correlation plot which illustrate various aspects of the color variable's dependence on the bar variable.  The software also generates 1-variable bar plots and 3-variable multi-panel plots. Note that the bar plot and box plot have fine scales along the bars that enable the determination of exact percentages, making tables unnecessary.  

A lot of thought and some highly nontrivial mathematics went into the design of these representations. For example, one may see from Figure~\ref{fig:2varPlot} that there appears to be a pattern in the bars. If the bar categories were presented in alphabetical order (as in Figure~\ref{fig:NARHS}, they would appear as a hodgepodge, with no discernible pattern.  The ordering of bars was determined by a ``traveling salesman'' algorithm, that minimizes the Manhattan distance (familiar  to mathematicians as the ``$\ell_1$ norm'') between successive bars, while secondarily maximizing the distance from the first to the last bar. 

Another problem that we had to solve was calculation speed. Since we are dealing with a database of 7 characteristics of over 30,000 students, to recalculate each graph from the original data was taking far too long. We required a system that could be used interactively, and ideally would generate graphs within a second or two. Fortunately we were able to solve this problem by reformatting the data into a much smaller multi-dimensional frequency table, on which all subsequent data operations were conducted. Data manipulations on this multidimensional table involved the operations of permuting and summing over different dimensions, meaning that some elementary combinatorics was also required.

It is true that these graphs are not fully self-explanatory, and  some explanations must be given in order for non-mathematical users to make sense of these pictures. But these explanations are not that complicated: even non-mathematical people have an idea of average, percentage, and range of values.  We are confident that additional development will bring additional simplifications.

\section{Conclusion}

We close with several suggestions to Nigerian academic mathematicians which may lead to a more robust environment for innovations with eventual practical impact. 
\begin{itemize}
\item
Resist the tendency to limit yourself to ``safe'' research within a narrow specialty. Within Nigerian mathematics there seems to be a very strong tendency towards narrow specialization. One key reason for this may be a dysfunctional system of academic promotion that rewards those who produce large numbers of  papers, regardless of the papers' depth or impact.  This situation is not limited to Nigeria, but is common to academic mathematics worldwide, including in the U.S., as William Stein has pointed out in  his presentations (see Figure~\ref{fig:steinSlide}).  

\begin{figure}[H]
	\centering
	\framebox{\includegraphics[width=3.5in]{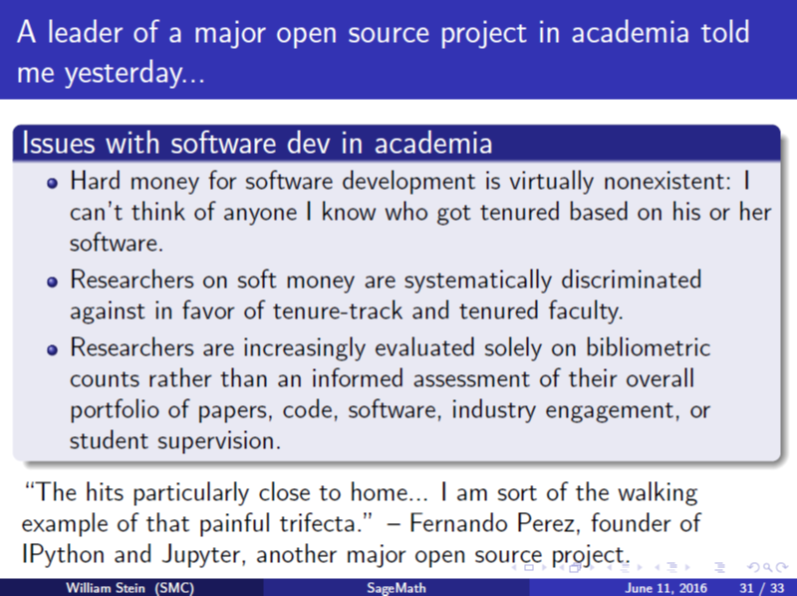}}
	\caption{Slide from William Stein's lecture at the University of Rochester (viewable on YouTube)\cite{Wing}, describing problems in the academic tenuring system.}
	\label{fig:steinSlide}
\end{figure}
\item
Resist the tendency towards mathematical elitism, namely the attitude that truly original mathematical research papers must be abstruse, complicated, and incomprehensible except to specialists. It is very easy for new mathematicians to try to make things unnecessarily complicated in order to make their work seem more "advanced". But this is exactly backwards: the purpose of mathematics is to simplify, not to obscure.
\item 
Develop fundamental intuitions, both in the curriculum and in research. This will assist in conveying mathematical ideas to people working in other areas.  For example, a great deal of linear algebra can be understood in terms of lines and planes in two or three dimensions.  Once  linear algebra is understood in these terms, then you can use these intuitions to explain Hilbert spaces, which is basically an infinite-dimensional real (or complex) vector space that also uses the Pythagorean theorem to define length and the inner product to define angle.  From Hilbert space you may move on to Banach spaces, which resemble Hilbert spaces but lack an inner product.  
\item
Recognize and promote mathematics department members who do productive work in high-impact areas, which may include (besides writing research papers) consulting, educating, and popularizing. Such members bring benefit to the entire department.  Encourage and reward those who are willing to take the risk not to follow ``conventional'' research. Be willing to accept dissertation proposals that are interdisciplinary and unconventional, particularly those that are oriented towards making an impact.
\item
Review the work of scientists, engineers and others that are working in math-related areas. Quite often they are using mathematical methods that they scarcely understand: and as a result they often either misuse methods, or employ them inefficiently. Or, they may be unaware of methods that you are quite familiar with.  As a mathematician, you have insights that others don't have. You can see coherence and structure  where others see a random jumble of equations. You can hear the music in the mathematics, whereas others hear only noise. 
\end{itemize}

\end{document}